\documentclass[12pt]{article}

\usepackage{amsmath}
\usepackage{amsthm,amscd,amsfonts}
\usepackage{amssymb, upref, color}
\usepackage{amsmath,amssymb,amsthm}
\usepackage{color}
\usepackage[dvips]{graphicx}
\usepackage{epsf,epsfig, verbatim} %subfigure,
\usepackage{latexsym,bm}
\usepackage{amsthm,amsmath,amssymb}
\usepackage{mathrsfs}
\usepackage{indentfirst}
\numberwithin{equation}{section}

\theoremstyle{plain}
\newtheorem{exam}{Example}[section]
\newtheorem{theorem}[exam]{Theorem}
\newtheorem{lemma}[exam]{Lemma}
\newtheorem{remark}[exam]{Remark}

\newtheorem{definition}[exam]{Definition}

\begin{document}
\date{}

%\newpage

\title{ %Linearization in non-autonomous differential equations based on algebraic  dichotomy
A  Hartman-Grobman theorem for algebraic dichotomies
\footnote{ This paper was jointly supported from the National Natural
Science Foundation of China under Grant (No. 11671176, 11931016), Grant Fondecyt (No. 1170466) and   Natural
Science Foundation of Zhejiang Province under Grant (No. LY20A010016).}
}
\author
{
Chaofan Pan$^{a}$\,\,\,\,\,
Manuel Pinto$^{b}$\,\, \,\,
Y.H. Xia$^{a}\footnote{Corresponding author. Y.H. Xia, xiaoutlook@163.com;yhxia@zjnu.cn. Address: College of Mathematics and Computer Science,  Zhejiang Normal University, 321004, Jinhua, China}$
\\
{\small \textit{$^a$ College of Mathematics and Computer Science,  Zhejiang Normal University, 321004, Jinhua, China}}\\
{\small \textit{$^b$ Departamento de Matem\'aticas, Universidad de Chile, Santiago, Chile }}\\
{\small Email: pancf0530@zjnu.edu.com;  pintoj.uchile@gmail.com; xiadoc@163.com; yhxia@zjnu.cn.}
}

\maketitle

\begin{abstract}
Algebraic dichotomy is a generalization of an exponential dichotomy (see Lin \cite{Lin}). This paper gives a version of Hartman-Grobman linearization theorem  assuming that linear system admits an algebraic dichotomy, which generalizes the Palmer's linearization theorem. Besides, we prove that the homeomorphism in the linearization theorem is  H\"{o}lder continuous (and has a H\"{o}lder continuous inverse). %To the best of our knowledge, the regularity of an equivalent function  appeared nowhere before in the published literature if linear system admits algebraic dichotomy.
Comparing with  exponential dichotomy, algebraic dichotomy is more complicate. {The exponential dichotomy leads us to the estimates $\int_{-\infty}^{t}e^{-\alpha(t-s)}ds$ and  $\int_{t}^{+\infty}e^{-\alpha(s-t)}ds$ which are convergent. However, the algebraic dichotomy will leads us to $\int_{-\infty}^{t}\left(\frac{\mu(t)}{\mu(s)}\right)^{-\alpha}ds$ or  $\int_{t}^{+\infty}\left(\frac{\mu(s)}{\mu(t)}\right)^{-\alpha}ds$, whose the convergence is unknown in the sense of Riemann.}



\end{abstract}
%Exponential dichotomy is a special case of algebraic dichotomy.
%\end{minipage}
%\end{center}
\section{Introduction}

\subsection{Brief history on dichotomy and $C^0$ linearization}

\indent  In 1930, Perron \cite{Perron} introduced the concept of the (classical or uniform) exponential dichotomy. Exponential dichotomy theory plays an important role in the differential equations. However, many scholars argued that  exponential dichotomy restricts many dynamic behaviors. For this reason, mathematicians proposed various concepts of the dichotomies which are more general than exponential dichotomy, for examples, ordinary dichotomy (see Coppel \cite{Coppel}), nonuniform exponential dichotomy (see Barreira and Valls \cite{Valls1,Valls2}, Barreira et al. \cite{Valls5}),  nonuniform polynomial dichotomy (see  Barreira and Valls \cite{Valls4}), $(h,k)$-dichotomy (see Naulin and Pinto \cite{Naulin},  Fenner and Pinto \cite{J. Fenner3}), nonuniform $(\mu,\nu)$-dichotomy (see Bento and Silva \cite{Silva}, Chang et al. \cite{Chang1}, Barreira et al. \cite{Valls3}, Chu \cite{Chu1}),  algebraic dichotomy (see Lin \cite{Lin}), $(h,k,\mu,\nu)$ dichotomy (see Zhang et al. \cite{J.M. Zhang, Zhang2}). In fact, the dichotomy was given for the different kinds of differential equations, such as the impulsive systems (see e.g \cite{Zhang2}), the difference equations (see e.g Barraria and Valls \cite{VallsD,VallsE}), the dynamic equations on time scales (see P\"otzche \cite{Potzche-JMAA,Potzche-NA}).  The most recently, the dichotomy theory was proposed and applied to the linear evolution equations with non-instantaneous impulsive effects (see Li et al. \cite{Li-Wang1, Li-Wang2}, Wang et al. \cite{Wang1,Wang-IMA}).
In this paper, we pay particular attention to the algebraic dichotomy introduced by
Lin \cite{Lin} discussed some basic properties of the algebraic dichotomy, and calculated the power of the weight function.\\
\indent On the other hand,  linear equations are mathematically well-understood but nonlinear systems are relatively difficult to investigate. For this reason, linearization of differential equations is very important. A basic contribution to the linearization problem for autonomous differential equations is the Hartman-Grobman theorem (see \cite{Hartman, Grobman}). Some improvements of the Hartman-Grobman theorem to infinite dimensional space can be found in Bates and Lu \cite{Bates1}, Hein and Pr\"{u}ss \cite{Hein-Pruss1}, Lu \cite{Lu1},
Pugh \cite{C.Pugh} and Reinfelds \cite{A.Reinfelds1,Reinfelds2}.
Palmer successfully generalized the Hartman-Grobman theorem to nonautonomous differential equations (see \cite{K.J.Palmer})
\begin{equation}
y'=A(t)y+f(t,y). \label{nueva_1}
\end{equation}
In order to weaken the conditions of Palmer's linearization theorem, some improvements were given in Backes et al. \cite{BDK-JDE} (without  exponential dichotomy), Barreira et al. \cite{L.Barreira2, L.Barreira3, L.Barreira4, L.Barreira9} (nonuniform dichotomy), Jiang \cite{L.Jiang1,L.Jiang2} (generalized dichotomy and ordinary dichotomy), %Reinfelds and Sermone \cite{ }, Sermone \cite{...}, Shi and Zhang \cite{Shi-Book},
Huerta \cite{Huerta2,Huerta1} (with nonuniform contraction), Papaschinopoulos \cite{Papaschinopoulos-A}, Pinto et al. \cite{Pinto-A}, Zou et al. \cite{Zou-Xia,Xia4}, Huang and Xia \cite{Xia5} (for the differential equations with piecewise constant argument), P\"otzche \cite{Potzche1} (for dynamic systems on time scales), Fenner and Pinto \cite{J.L.Fenner2}, Xia et al. \cite{Xia1,Xia3} and Zhang et al. \cite{J.M. Zhang,Zhang2} (for the instantaneous impulsive system). %\color {red}{(nonuniform $``(h,k,\mu,\nu)"$ dichotomy for ordinary differential equations and nonautonomous impulsive differential equations)}.

\subsection{Motivation and novelty}

In this paper,  we  pay particular attention to the effect of the algebraic dichotomy imposing on the linearization of the differential equations.
Palmer's linearization theorem requires two essential conditions: (i) the nonlinear term $f$ is uniformly bounded and Lipschitzian; (ii) the linear system
\begin{equation}
 x'(t)=A(t)x(t) \label{eq1}
\end{equation}
possesses an exponential dichotomy.
In the present paper, we try to reduce the second condition. Motivated by Lin's algebraic dichotomy (see Lin \cite{Lin}) and the works of Palmer \cite{K.J.Palmer} and Zhang et al \cite{J.M. Zhang}, we study the $C^0$ linearization with the algebraic dichotomy.  Further more, we prove that  the homeomorphism and its inverse are H\"older continuous under the assumption of the algebraic dichotomy.  When the algebraic dichotomy reduces to the exponential dichotomy, our results generalize and improve the previous ones.
 Comparing with  exponential dichotomy, algebraic dichotomy is more general. {The exponential dichotomy leads to the estimates $\int_{-\infty}^{t}e^{-\alpha(t-s)}ds$ and  $\int_{t}^{+\infty}e^{-\alpha(s-t)}ds$ which are convergent. However, the algebraic dichotomy will leads us to $\int_{-\infty}^{t}\left(\frac{\mu(t)}{\mu(s)}\right)^{-\alpha}ds$ or  $\int_{t}^{+\infty}\left(\frac{\mu(s)}{\mu(t)}\right)^{-\alpha}ds$, whose the convergence is unknown in the sense of Riemann.} This brings more difficulties to our research.

\indent  The structure of our paper  as follows.  In Section 2, we give our main results. In Section 3, we give some preliminary results. In Section 4, we give  rigorous proofs to show the regularity of the equivalent function $H(t, x)$ and $G(t, y)$. Finally, we give an example to illustrate our linearization theorem.

\section{Statement of main results}
Consider the following two non-autonomous systems
\begin{equation}\label{2.1}
x'=A(t)x+f(t,x)
\end{equation}
and
\begin{equation}\label{2.2}
x'=A(t)x,
\end{equation}
where $x\in \mathbb{R}^{n}$, $t\in \mathbb{R}$. Suppose that $A(t)$ is a $n\times n$ continuous and bounded matrix defined on $\mathbb{R}$. Let $T(t, s)$ be the evolution operator of system (\ref{2.2}) satisfying $T(t,s)x(s)=x(t)$, $t,s\in\mathbb{R}$, for any solution $x(t)$ of system (\ref{2.2}).  Clearly, $T(t,t)=Id$ and
\begin{align}\label{2.4}
T(t,\tau)T(\tau,s)=T(t,s), \ \ \ t,s,\tau\in\mathbb{R}.
\end{align}
An increasing function $\mu:\mathbb{R}\to\mathbb{R}_{0}^{+}$ is said to be a growth rate if $\mu(0)=1$,
\begin{align}\label{for2.5}
\lim_{t \to -\infty}\mu(t)=0 \quad \mathrm{and} \quad \lim_{t \to +\infty}\mu(t)=+\infty.
\end{align}
In the following, we always  assume that $\mu(t)$ is growth rate.
\begin{definition}\label{def2.2}
\cite{Lin} Linear system (\ref{2.2}) is said to admit an algebraic dichotomy, if there exists a projection $P(s)$ and constants $K > 0$, $\alpha > 0$ such that
\begin{equation}\label{2.6}
\begin{split}
T (t, s)P(s) &=P(t)T (t, s),
\end{split}
\end{equation}
\begin{equation}\label{2.7}
\begin{split}
||T(t,s)P(s)||&\leq K\left(\frac{\mu(t)}{\mu(s)}\right)^{-\alpha},
\quad \mathrm{if} \quad t\geq s,\\
||T(t,s)Q(s)||&\leq K\left(\frac{\mu(s)}{\mu(t)}\right)^{-\alpha},
\quad \mathrm{if} \quad t\leq s
\end{split}
\end{equation}
hold, where $P(s)+Q(s)=I$.
\end{definition}
Obviously, $\mu(t)=e^{t}$ yields an exponential dichotomy.  %Therefore, algebraic dichotomy systems are certainly more general than exponential dichotomic systems.

To show its generality, we give an example of algebraic dichotomy.
\begin{exam}\label{exam2.2}
\rm Consider the differential equation in $\mathbb{R}^2$.
\begin{equation}\label{exam2.8}
\begin{split}
\begin{cases}
x'_{1}=-\eta_{1}\left(\frac{\mu'(t)}{\mu(t)}\right)x_{1},\\
x'_{2}=\eta_{2}\left(\frac{\mu'(t)}{\mu(t)}\right)x_{2},
\end{cases}
\end{split}
\end{equation}
\rm for $t\in\mathbb{R}$, where $\eta_{1}, \eta_{2}$ are positive constants.
\end{exam}
\noindent From system (\ref{exam2.8}), we get
\begin{align*}
x_{1}=\left(\frac{\mu'(t)}{\mu(t)}\right)^{-\eta_{1}}V_{1},\quad
x_{2}=\left(\frac{\mu'(t)}{\mu(t)}\right)^{\eta_{2}}V_{2},
\end{align*}
where $V_{1}, V_{2}$ are constants.\\
Taking $P(s)=diag\{1,0\}$, $Q(s)=diag\{0,1\}$.
Then, we have
\begin{equation}
T(t,s)P(s)=
\left(
\begin{array}{cc}
\left(\frac{\mu(t)}{\mu(s)}\right)^{-\eta_{1}}&0\\
0&0\\
\end{array}
\right),\quad \mathrm{for} \quad t\geq s,
\end{equation}
\begin{equation}
T(t,s)Q(s)=
\left(
\begin{array}{cc}
0&0\\
0&\left(\frac{\mu(s)}{\mu(t)}\right)^{-\eta_{2}}\\
\end{array}
\right),\quad \mathrm{for} \quad t\leq s.
\end{equation}
Taking $\eta_{3}=\max\{\eta_{1},\eta_{2}\}$. It follows that
\begin{align*}
||T(t,s)P(s)||\leq\left(\frac{\mu(t)}{\mu(s)}\right)^{-\eta_{3}},\quad \mathrm{for} \quad t\geq s,
\end{align*}
\begin{align*}
||T(t,s)Q(s)||\leq\left(\frac{\mu(s)}{\mu(t)}\right)^{-\eta_{3}},\quad \mathrm{for} \quad t\leq s.
\end{align*}
This implies that equation (\ref{exam2.8}) admits an algebraic dichotomy with
\begin{align*}
K=1,\quad \alpha=\eta_{3}.
\end{align*}
\begin{remark}
Particularly, if $\mu(t)=e^{t}$, Example \ref{exam2.2} implies that equation (\ref{exam2.8}) admits an exponential dichotomy with $K=1, \alpha=\eta_{3}$.
\end{remark}
\begin{definition}
\cite{J. Fenner3} Linear system (\ref{2.2}) is said to admit an $(h,k)$ dichotomy, if there exists  a projection $P(s)$, constants $K > 0$, $\alpha > 0$ and piecewise right continuous function $h(t)$, $k(t)$  such that
\begin{equation}
\begin{split}
\begin{cases}
||U(t,s)P(s)||\leq Kh(t)h^{-1}(s)e^{-\alpha(t-s)},\quad t\geq s,\\
||U(t,s)Q(s)||\leq Kk(t)k^{-1}(s)e^{-\alpha(s-t)},\quad t\leq s.
\end{cases}
\end{split}
\end{equation}
hold, where $P(s)+Q(s)=I$.
\end{definition}
If $h(t)=k(t)=\mu(t)$, then we take $\max\{e^{-\alpha(t-s)},t\geq s\}$ and $\max\{e^{-\alpha(t-s)},t\leq s\}$, which are 1.
Thus, we get an algebraic dichotomy.
\begin{definition}
\cite{J. Fenner3} We say that $h$ and $k$ are fulfill a compensation law on $\mathbb{R}$ if there exists a positive constant $C_{h,k}$ such that
\begin{equation*}
k(t)k^{-1}(s)\leq C_{h,k}h(t)h^{-1}(s),\quad t\geq s.
\end{equation*}
The system is said to be a $h$-system, if it has a $(h, h)$ dichotomy.
\end{definition}
Clearly, a system having an $(h,k)$-dichotomy with compensation law belongs to the
class of $h$-systems.
\begin{definition}\label{def2.1}
Suppose that there exists a function $H : \mathbb{R} \times \mathbb{R}^{n} \to \mathbb{R}^{n} $such that\\
$(i)$ for each fixed $t$, $H(t,\cdot)$ is a homeomorphism of $\mathbb{R}^{n}$ into $\mathbb{R}^{n}$;\\
$(ii)$ $||H(t, x)-x||$ is uniformly bounded with respect to $t$;\\
$(iii)$ assume that $G(t, \cdot) = H^{-1}(t, \cdot)$ also has property $(ii)$;\\
$(iv)$ if $x(t)$ is a solution of system (\ref{2.1}), then $H(t, x(t))$ is a solution of system (\ref{2.2}); and if $y(t)$ is a solution of system (\ref{2.2}), then $G(t, y(t))$ is a solution
of system (\ref{2.1}).\\
\indent If such a map $H_{t}(H_{t} := H(t, \cdot))$ exists, then  system (\ref{2.1}) is topologically conjugated to  system (\ref{2.2}) and the transformation $H(t, x)$ is called an equivalent function.
\end{definition}
Now we are in a position to state our main results.

\begin{theorem}\label{thm3.1}
Suppose that system (\ref{2.2}) admits an algebraic dichotomy, and $f(t,x)$ satisfies
\begin{align}\label{3.1}
||f(t,x)|| &\leq \beta \mu'(t)\mu^{-1}(t),\\
||f(t,x_{1})-f(t,x_{2})|| &\leq \gamma \mu'(t)\mu^{-1}(t)||x_{1}-x_{2}||,\\
6K\gamma\alpha^{-1}&<1.
\end{align}
Then  nonlinear system (\ref{2.1}) is topologically conjugated to their linear part $\dot{x}=A(t)x$, and the equivalent function
\begin{align*}
||H(t,x)-x||\leq2K\beta\alpha^{-1}.
\end{align*}
Denote $H^{-1}(t,\cdot)=G(t,\cdot)$, then $G(t,y)$ also satisfies
\begin{align*}
||G(t,y)-y||\leq2K\beta\alpha^{-1}.
\end{align*}
\end{theorem}
\begin{theorem}\label{thm3.2}
Suppose that the conditions in Theorem \ref{thm3.1} are satisfied. Moreover, assuming that %$||\mu||_{0}=\sup_{t\in\mathbb{R}}\int_{t}^{t+1}\mu'(s)\mu^{-1}(s)ds$
$\alpha>\gamma$.
Then there exist constants $p,q >0$, $0<p',q'<1$ such that
\begin{align}
||H(t,x)-H(t,x')||\leq p||x-x'||^{q}, \quad \mathrm{if} \ ||x-x'||<1,\\ ||G(t,y)-G(t,y')||\leq p'||y-y'||^{q'}, \quad  \mathrm{if} \ ||y-y'||<1,
\end{align}
where $G(t, \cdot) = H^{-1}(t, \cdot)$.
\end{theorem}
\begin{remark}
Assuming that $\mu(t)=e^{t}$ in Theorem \ref{thm3.1}, Theorem \ref{thm3.1} reduces to the classical Palmer linearization theorem, (see \cite{K.J.Palmer}). However, Palmer did not study the regularity  of equivalent function $H(t, x)$. We remark that there are good results for the linearization of $(h,k,\mu,\nu)$-dichotomy (\cite{J.M. Zhang,Zhang2}). However, they did not study the regularity of homeomorphisms mapping the nonlinear systems onto its linearization.
\end{remark}
%\begin{corollary}
%Suppose that system (\ref{2.2}) admits a $(h,k)$ dichotomy, and $f(t,x)$ satisfies
%\begin{align}\label{3.1}
%||f(t,x)|| &\leq \beta \min\{h'(t)h^{-1}(t),k'(t)k^{-1}(t)\},\\
%||f(t,x_{1})-f(t,x_{2})|| &\leq \gamma \min\{h'(t)h^{-1}(t),k'(t)k^{-1}(t)\}||x_{1}-x_{2}||,\\
%6K\gamma\alpha^{-1}&<1.
%\end{align}
%Then  nonlinear system (\ref{2.1}) is topologically conjugated to their linear part $\dot{x}=A(t)x$.
%\end{corollary}

\section{Preliminary results}
We split the proof of Theorem \ref{thm3.1} into several lemmas. Suppose that $X(t,t_{0},x_{0})$ is the solution of  system (\ref{2.1}) with the initial value condition $X(t_{0})=x_{0}$,
and  $Y(t,t_{0},y_{0})$ is the solution of  system (\ref{2.2}) with the initial value condition $Y(t_{0})=y_{0}$.

We start with a fundamental lemma, which shows that linear system (\ref{2.2}) has no other bounded solutions except for the zero solution under our hypothesis.
\begin{lemma}\label{lem4.1}
 If  linear system (\ref{2.2}) has an algebraic dichotomy, then the bounded solution of linear system (\ref{2.2}) is the zero solution.
 \end{lemma}
{\bf Proof}.
Let $x(t)$ be the bounded solution of system (\ref{2.2}). There exists $n$ order real vector $x(s)$, such that
\begin{align}\label{cy4.1}
x(t)=T(t,s)x(s).
\end{align}
From (\ref{cy4.1}), we get
\begin{align*}
x(t)=T(t,s)P(s)x(s)+T(t,s)Q(s)x(s).
\end{align*}
If $P(s)x(s)\neq0$, consider the case  $t\leq s$. We have
\begin{align*}
||x(t)||=||T(t,s)P(s)x(s)+T(t,s)Q(s)x(s)||\geq||T(t,s)P(s)x(s)||-||T(t,s)Q(s)x(s)||.
\end{align*}
From (\ref{2.4}), (\ref{2.6}) and (\ref{2.7}), we get
\begin{align*}
||P(s)x(s)||&=||T(s,s)P(s)T(s,t)T(t,s)P(s)x(s)||\nonumber\\
&\leq||T(s,t)P(t)|| \ ||T(t,s)P(s)x(s)||\nonumber\\
&\leq K\left(\frac{\mu(s)}{\mu(t)}\right)^{-\alpha}\ ||T(t,s)P(s)x(s)||.
\end{align*}
Thus, we obtain
\begin{align*}
||T(t,s)P(s)x(s)||\geq\frac{1}{K}\left(\frac{\mu(s)}{\mu(t)}\right)^{\alpha}||P(s)x(s)||.
\end{align*}
From (\ref{2.7}), we get
\begin{align*}
||T(t,s)Q(s)x(s)||\leq K\left(\frac{\mu(s)}{\mu(t)}\right)^{-\alpha}||x(s)||.
\end{align*}
Hence,
\begin{align*}
||x(t)|| \geq\frac{1}{K}\left(\frac{\mu(s)}{\mu(t)}\right)^{\alpha}||P(s)x(s)||-
K\left(\frac{\mu(s)}{\mu(t)}\right)^{-\alpha}||x(s)||.
\end{align*}
From (\ref{for2.5}), when $t \to -\infty$, $\mu(t) \to 0$. Thus,
\begin{align*}
\lim_{t \to -\infty}||x(t)||=+\infty.
\end{align*}
This contradicts with the boundedness of $x(t)$, thus $P(s)x(s)=0$. Similarly, we  get $Q(s)x(s)=0$, when $t>s$. Thus, $x(t)=0$, $t\in \mathbb{R}$.\\
\indent This conclusion is key to the following lemmas.
\begin{lemma}\label{lem4.2}
For any given $(\tau,\xi)$, system
\begin{equation}\label{equ4.1}
Z'(t)=A(t)Z(t)-f(t,X(t,\tau,\xi))
\end{equation}
has a unique bounded solution $h(t,(\tau,\xi))$, and $h(t,(\tau,\xi))\leq2 K\beta\alpha^{-1}$.
\end{lemma}
{\bf Proof.} For any given $(\tau,\xi)$, taking
\begin{align*}
Z_{0}(t)=-\int_{-\infty}^{t}T(t,s)P(s)f(s,X(s,\tau,\xi))ds
+\int_{t}^{+\infty}T(t,s)Q(s)f(s,X(s,\tau,\xi))ds.
\end{align*}
Differentiating $Z_{0}(t)$, we get $Z_{0}(t)$ is a solution of system (\ref{equ4.1}).\\
\indent Now we prove that $Z_{0}(t)$ is the unique bounded solution of system (\ref{equ4.1}).
\begin{align*}
||Z_{0}(t)||& \leq \int_{-\infty}^{t}K\left(\frac{\mu(t)}{\mu(s)}\right)^{-\alpha}\beta\mu'(s)\mu^{-1}(s)ds
+\int_{t}^{+\infty}K\left(\frac{\mu(s)}{\mu(t)}\right)^{-\alpha}\beta\mu'(s)\mu^{-1}(s)ds \nonumber\\
&\leq K\beta\mu^{-1}(t)\int_{-\infty}^{t}\mu^{\alpha-1}(s)\mu'(s)ds+K\beta\mu^{\alpha}(t)
\int_{t}^{+\infty}\mu^{-\alpha-1}(s)\mu'(s)ds \nonumber \\
& \leq 2 K\beta \alpha^{-1}.
\end{align*}
For any given $(\tau,\xi)$, system (\ref{equ4.1}) is a linear non-homogeneous system, and its linear part
\begin{equation*}
Z'(t)=A(t)Z(t)
\end{equation*}
has an algebraic dichotomy. According to Lemma \ref{lem4.1},
system (\ref{equ4.1}) has a unique bounded solution $Z_{0}(t)$.
Since $Z_{0}(t)$ is related to $\tau,\xi$,  we denote $Z_{0}(t)$ as $h(t,(\tau,\xi))$. From the above proof, we get
$h(t,(\tau,\xi))\leq2K\beta\alpha^{-1}$.\\

\begin{lemma}\label{lem4.3}
For any given $(\tau,\xi)$, system
\begin{equation}\label{equ4.2}
Z'(t)=A(t)Z(t)+f(t,Y(t,\tau,\xi)+Z(t))
\end{equation}
has a unique bounded solution $g(t,(\tau,\xi))$, and $g(t,(\tau,\xi))\leq2K\beta\alpha^{-1}$.
\end{lemma}
{\bf Proof.} Let
\begin{align*}
\Omega:=\{z:\mathbb{R}\to X\ |\ ||z|| \leq 2K\beta\alpha^{-1}\}.
\end{align*}
Defining the following mapping:
\begin{align*}
\mathscr{F} z=\int_{-\infty}^{t}T(t,s)P(s)f(s,Y(s,\tau,\xi)+z)ds
-\int_{t}^{+\infty}T(t,s)Q(s)f(s,Y(s,\tau,\xi)+z)ds.
\end{align*}
Furthermore,
\begin{align*}
||\mathscr{F} z||\leq \int_{-\infty}^{t}K\left(\frac{\mu(t)}{\mu(s)}\right)^{-\alpha}\beta\mu'(s)\mu^{-1}(s)ds
+\int_{t}^{+\infty}K\left(\frac{\mu(s)}{\mu(t)}\right)^{-\alpha}\beta\mu'(s)\mu^{-1}(s)ds
\leq2K\beta\alpha^{-1}.
\end{align*}
Thus, $\mathscr{F}$ is a self-mapping in $\Omega$. Moreover,
\begin{align*}
||\mathscr{F} z_{1}-\mathscr{F} z_{2}||\leq& \int_{-\infty}^{t}K\left(\frac{\mu(t)}{\mu(s)}\right)^{-\alpha}\gamma\mu'(s)\mu^{-1}(s)||z_{1}-z_{2}||ds\\
&+\int_{t}^{+\infty}K\left(\frac{\mu(s)}{\mu(t)}\right)^{-\alpha}\gamma\mu'(s)\mu^{-1}(s)||z_{1}-z_{2}||ds
\nonumber\\
\leq&2K\gamma\alpha^{-1}||z_{1}-z_{2}||\nonumber\\
\leq&\frac{1}{3}||z_{1}-z_{2}||.
\end{align*}
\indent From Banach's fixed point theorem, $\mathscr{F}$ has a unique fixed point $Z_{1}(t)$ in $\Omega$ , and $Z_{1}(t)$ satisfies
\begin{align*}
Z_{1}(t)=\int_{-\infty}^{t}T(t,s)P(s)f(s,Y(s,\tau,\xi)+Z_{1}(s))ds
-\int_{t}^{+\infty}T(t,s)Q(s)f(s,Y(s,\tau,\xi)+Z_{1}(s))ds.
\end{align*}
Differentiating $Z_{1}(t)$, we get $Z_{1}(t)$ is a solution of equation (\ref{equ4.2}). Since
$||Z_{1}(t)||\leq 2K\beta\alpha^{-1}$,  $Z_{1}(t)$ is a bounded solution of equation (\ref{equ4.2}).\\
\indent Now we prove the uniqueness. Let $Z_{2}(t)$ be another bounded solution. By the variation formula, we get
\begin{align*}
Z_{2}(t)=T(t,0)x_{0}&+\int_{0}^{t}T(t,s)f(s,Y(s,\tau,\xi)+Z_{2}(s))ds\\
=T(t,0)x_{0}&+\int_{0}^{t}T(t,s)P(s)f(s,Y(s,\tau,\xi)+Z_{2}(s))ds\\
&+\int_{0}^{t}T(t,s)Q(s)f(s,Y(s,\tau,\xi)+Z_{2}(s))ds\\
=T(t,0)x_{0}&+\int_{-\infty}^{t}T(t,s)P(s)f(s,Y(s,\tau,\xi)+Z_{2}(s))ds\\
&-\int_{-\infty}^{0}T(t,s)P(s)f(s,Y(s,\tau,\xi)+Z_{2}(s))ds\\
&+\int_{0}^{+\infty}T(t,s)Q(s)f(s,Y(s,\tau,\xi)+Z_{2}(s))ds\\
&-\int_{t}^{+\infty}T(t,s)Q(s)f(s,Y(s,\tau,\xi)+Z_{2}(s))ds.
\end{align*}
Since
\begin{align*}
& \int_{-\infty}^{0}T(t,s)P(s)f(s,Y(s,\tau,\xi)+Z_{2}(s))ds\\
=&T(t,0)T(0,t)\int_{-\infty}^{0}T(t,s)P(s)f(s,Y(s,\tau,\xi)+Z_{2}(s))ds,
\end{align*}
and
\begin{equation}\label{cy4.4}
\begin{split}
&||T(0,t)\int_{-\infty}^{0}T(t,s)P(s)f(s,Y(s,\tau,\xi)+Z_{2}(s))ds||\\
=&||\int_{-\infty}^{0}T(0,t)T(t,s)P(s) f(s,Y(s,\tau,\xi)+Z_{2}(s))ds||\\
=&||\int_{-\infty}^{0}T(0,s)P(s)f(s,Y(s,\tau,\xi)+Z_{2}(s))ds||\\
\leq& \int_{-\infty}^{0}K\left(\frac{\mu(0)}{\mu(s)}\right)^{-\alpha}\beta\mu'(s)\mu^{-1}(s)ds\\
\leq&K\beta\alpha^{-1}.
\end{split}
\end{equation}
From (\ref{cy4.4}), we get $T(0,t)\int_{-\infty}^{0}T(t,s)P(s)f(s,Y(s,\tau,\xi)+Z_{2}(s))ds$ is convergent, and denoting it as $x_{1}$. Thus, we get
\begin{align*}
\int_{-\infty}^{0}T(t,s)P(s)f(s,Y(s,\tau,\xi)+Z_{2}(s))ds=T(t,0)x_{1}.
\end{align*}
Similarly,
\begin{align*}
\int_{0}^{+\infty}T(t,s)Q(s)f(s,Y(s,\tau,\xi)+Z_{2}(s))ds=T(t,0)x_{2}.
\end{align*}
Therefore,
\begin{align*}
Z_{2}(t)=T(t,0)(x_{0}+x_{1}+x_{2})&+\int_{-\infty}^{t}T(t,s)P(s)f(s,Y(s,\tau,\xi)+Z_{2}(s))ds\\
&-\int_{t}^{+\infty}T(t,s)Q(s)f(s,Y(s,\tau,\xi)+Z_{2}(s))ds.
\end{align*}
Since $Z_{2}(t)$ is a bounded solution, and
\begin{align*}
 \int_{-\infty}^{t}T(t,s)P(s)f(s,Y(s,\tau,\xi)+Z_{2}(s))ds
-\int_{t}^{+\infty}T(t,s)Q(s)f(s,Y(s,\tau,\xi)+Z_{2}(s))ds
\end{align*}is bounded.
Therefore, $T(t,0)(x_{0}+x_{1}+x_{2})$ is also bounded.\\
\indent Since $T(t,0)(x_{0}+x_{1}+x_{2})$ is the bounded solution of $Z'=A(t)Z$. From Lemma \ref{lem4.1},  we get
 $T(t,0)(x_{0}+x_{1}+x_{2})=0$. Thus,
\begin{align*}
Z_{2}(t)=\int_{-\infty}^{t}T(t,s)P(s)f(s,Y(s,\tau,\xi)+Z_{2}(s))ds
-\int_{t}^{+\infty}T(t,s)Q(s)f(s,Y(s,\tau,\xi)+Z_{2}(s))ds.
\end{align*}
Moreover,
\begin{align*}
||Z_{1}(t)-Z_{2}(t)||\leq&\int_{-\infty}^{t}||T(t,s)P(s)||\beta \mu'(s) \mu ^{-1}(s)||Z_{1}(s)-Z_{2}(s)||ds\nonumber\\
&+\int_{t}^{+\infty}||T(t,s)Q(s)||\beta \mu'(s) \mu ^{-1}(s)||Z_{1}(s)-Z_{2}(s)||ds\nonumber\\
\leq&2K\gamma\alpha^{-1}||Z_{1}(t)-Z_{2}(t)||\nonumber\\
\leq&\frac{1}{3}||Z_{1}(t)-Z_{2}(t)||.
\end{align*}
Thus, $Z_{1}(t)=Z_{2}(t)$. The bounded solution of system (\ref{equ4.2}) is unique. This solution is related to $(\tau,\xi)$, denoting it as $g(t,(\tau,\xi))$ .
From the above proof, we get
\begin{align*}
g(t,(\tau,\xi))\leq2K\beta\alpha^{-1}.
\end{align*}

\begin{lemma}\label{lem4.4}
Let $x(t)$ be any solution of  system (\ref{2.1}), system
\begin{equation}\label{equ4.3}
Z'(t)=A(t)Z(t)+f(t,x(t)+Z(t))-f(t,x(t))
\end{equation}
has a unique bounded solution $Z=0$.
\end{lemma}
{\bf Proof.} Obviously, $Z=0$ is a bounded solution of  system (\ref{equ4.3}). Next, we prove the uniqueness of the bounded solution. Let $Z_{3}(t)$ be another bounded solution. By the variation formula, we get
\begin{align*}
Z_{3}(t)=T(t,0)x(0)+\int_{0}^{t}T(t,s)[f(s,x(s)+Z_{3}(s))-f(s,x(s))]ds.
\end{align*}
Similar to the proof of Lemma \ref{lem4.3}, we get
\begin{align*}
Z_{3}(t)=&\int_{-\infty}^{t}T(t,s)P(s)[f(s,x(t))+Z_{3}(s))-f(s,x(s))]ds\\
&-\int_{t}^{+\infty}T(t,s)Q(s)[f(s,x(t))+Z_{3}(s))-f(s,x(s))]ds.
\end{align*}
Moreover,
\begin{align*}
||Z_{3}(t)||\leq &\int_{-\infty}^{t}||T(t,s)P(s)||\beta \mu'(s) \mu ^{-1}(s)||Z_{3}(s)||ds\\
&+\int_{t}^{+\infty}||T(t,s)Q(s)||\beta \mu'(s) \mu ^{-1}(s)||Z_{3}(s)||ds\nonumber\\
\leq &2K\gamma\alpha^{-1}||Z_{3}(t)||\nonumber\\
\leq&\frac{1}{3}||Z_{3}(t)||.
\end{align*}
Thus, $Z_{3}(t)=0$.

Now we construct two functions as follows:
\begin{align}
H(t,x)=x+h(t,(t,x)) \label{4.26},\\
G(t,y)=y+g(t,(t,y))\label{4.27}.
\end{align}

\begin{lemma}\label{lem4.5}
For any given $(t_{0},x_{0})$, $H(t,X(t,t_{0},x_{0}))$ is the solution of  linear system (\ref{2.2}).
\end{lemma}
{\bf Proof.}
Replacing $(\tau,\xi)$ in system (\ref{equ4.1}) with $(t,X(t,\tau,\xi))$. From  the uniqueness of the bounded solution, we obtain
\begin{align}\label{4.6}
H(t,X(t,t_{0},x_{0}))=X(t,t_{0},x_{0})+h(t,(t_{0},x_{0})).
\end{align}
Differentiating (\ref{4.6}), we get
\begin{align*}
[H(t,X(t,t_{0},x_{0}))]'&=A(t)X(t,t_{0},x_{0})+f(t,X(t,t_{0},x_{0}))
+A(t)h(t,(t_{0},x_{0}))-f(t,X(t,t_{0},x_{0}))\nonumber\\
&=A(t)H(t,X(t,t_{0},x_{0})).
\end{align*}
This shows that $H(t,X(t,t_{0},x_{0}))$ is the solution of  linear system (\ref{2.2}).

\begin{lemma}\label{lem4.6}
For any given $(t_{0},y_{0})$, $G(t,Y(t,t_{0},y_{0}))$ is the solution of  system (\ref{2.1}).
\end{lemma}
{\bf Proof.} Replacing  $(\tau,\xi)$ in system (\ref{equ4.2}) with $(t,Y(t,\tau,\xi))$. From Lemma \ref{lem4.3} and Lemma \ref{4.27}, we get
\begin{align}\label{4.7}
G(t,Y(t,t_{0},y_{0}))=Y(t,t_{0},y_{0})+g(t,(t_{0},y_{0})).
\end{align}
Differentiating (\ref{4.7}), we get
\begin{align*}
[G(t,Y(t,t_{0},y_{0}))]'&=A(t)Y(t,t_{0},y_{0})
+A(t)g(t,(t_{0},g_{0}))+f(t,Y(t,t_{0},y_{0})+g(t,(t_{0},g_{0})))\nonumber\\
&=A(t)G(t,Y(t,t_{0},y_{0}))+f(t,G(t,Y(t,t_{0},y_{0})))
\end{align*}
This shows that $G(t,Y(t,t_{0},y_{0}))$ is the solution of system (\ref{2.1}).

\begin{lemma}\label{lem4.7}
For any $t\in \mathbb{R}$, $y\in \mathbb{R}^{n}$, we has always
\begin{align}
H(t,G(t,y))=y.
\end{align}
\end{lemma}
{\bf Proof.}
Let $y(t)$ be any solution of system (\ref{2.2}). From Lemma \ref{lem4.6}, we know that $G(t,y(t))$ is the solution of system (\ref{2.1}). From Lemma \ref {lem4.5}, we know that $H(t,G(t,y(t)))$ is the solution of system (\ref{2.2}), denoting it as $y_{0}(t)$.
Let $M(t)=y(t)-y_{0}(t)$. Then, we have
$$M'(t)=y'(t)-y'_{0}(t)=A(t)y(t)-A(t)y_{0}(t)=A(t)M(t).$$
Thus, $M(t)$ is the solution of system (\ref{2.2}).
Moreover,
\begin{align*}
||M(t)||&=||y(t)-y_{0}(t)||\\
&=||y(t)-H(t,G(t,y(t)))||\nonumber\\
&\leq||y(t)-G(t,y(t))||+||G(t,y(t))-H(t,G(t,y(t)))||\nonumber\\
&=||g(t,(t,y(t)))||+||h(t,(t,G(t,y(t))))||\nonumber\\
&\leq 2K\beta\alpha^{-1}+ 2K\beta\alpha^{-1}\nonumber\\
&=4K\beta\alpha^{-1}.
\end{align*}
This shows that $M(t)$ is a bounded solution of $x'=A(t)x$. From Lemma \ref{lem4.1}, we get $M(t)=0$.
Thus,
\begin{align*}
y(t)=y_{0}(t),\ \ \ \ H(t,G(t,y))=y.
\end{align*}

\begin{lemma}\label{lem4.8}
For any $t\in \mathbb{R},x\in \mathbb{R}^{n}$, we has always
\begin{align}
G(t,H(t,x))=x.
\end{align}
\end{lemma}
{\bf Proof.}
Let $x(t)$ be any solution of system(\ref{2.1}). From Lemma \ref{lem4.5},  we know that $H(t,x(t))$ is the solution of system (\ref{2.2}). From Lemma \ref{lem4.6}, we get $G(t,H(t,x(t)))$ is the solution of system(\ref{2.1}). Denoting it as $x_{0}(t)$.\\
\indent Let $N(t)=x_{0}(t)-x(t)$, differentiating it, we get
\begin{align*}
N'(t)&=x_{0}'(t)-x'(t)\\
&=A(t)x_{0}(t)+f(t,x_{0}(t))-A(t)x(t)-f(t,x(t))\\
&=A(t)N(t)+f(t,x(t)+N(t))-f(t,x(t)).
\end{align*}
Thus, $N(t)$ is a solution of  system (\ref{equ4.3}). Moreover,
\begin{align*}
||N(t)||&=||x(t)-x_{0}(t)||\\
&=||x(t)-G(t,H(t,x(t)))||\nonumber\\
&\leq||x(t)-H(t,x(t))||+||H(t,x(t))-G(t,H(t,x(t)))||\nonumber\\
&=||h(t,(t,x(t)))||+||g(t,(t,h(t,x(t))))||\nonumber\\
&\leq 2K\beta\alpha^{-1}+ 2K\beta\alpha^{-1}\nonumber\\
&=4K\beta\alpha^{-1}.
\end{align*}
This shows that $N(t)$ is a bounded solution of  system (\ref{equ4.3}). From Lemma \ref{lem4.4}, we have
 $N(t)=0$. Thus,
\begin{align*}
x(t)=x_{0}(t),\ \ \ \ G(t,H(t,x))=x.
\end{align*}
\begin{lemma}\label{lem4.10}
Denoting that $\sup\limits_{t\in \mathbb{R}}||A(t)||=M$. Then we have
\begin{align}
||X(t,t_{0},x_{0})-X(t,t_{0},x'_{0})||&\leq||x_{0}-x'_{0}||e^{M(t-t_{0})}
\left(\frac{\mu(t)}{\mu(t_{0})}\right)^{\gamma},\\
||Y(t,t_{0},y_{0})-Y(t,t_{0},y'_{0})||&\leq||y_{0}-y'_{0}||e^{M(t-t_{0})}.
\end{align}
\end{lemma}
{\bf Proof.} By the variation formula, we get
\begin{align*}
X(t,t_{0},x_{0})&=x_{0}+\int_{t_{0}}^{t}A(s)X(s,t_{0},x_{0})+f(s,X(s,t_{0},x_{0}))ds,\\
Y(t,t_{0},y_{0})&=y_{0}+\int_{t_{0}}^{t}A(s)Y(s,t_{0},x_{0})ds.
\end{align*}
Then,
\begin{align*}
||X(t,t_{0},x_{0})-X(t,t_{0},x'_{0})||&\leq||x_{0}-x'_{0}||+
\int_{t_{0}}^{t}M||X(s,t_{0},x_{0})-X(s,t_{0},x'_{0})||\\&+
\gamma
\mu'(s)\mu^{-1}(s)||X(s,t_{0},x_{0})-X(s,t_{0},x'_{0})||ds,\\
||Y(t,t_{0},y_{0})-Y(t,t_{0},y'_{0})||&\leq||y_{0}-y'_{0}||+
\int_{t_{0}}^{t}M||Y(s,t_{0},y_{0})-Y(s,t_{0},y'_{0})||ds.
\end{align*}
From Bellman's inequality, we get
\begin{align*}
||X(t,t_{0},x_{0})-X(t,t_{0},x'_{0})||&\leq ||x_{0}-x'_{0}|| e^{\int_{t_{0}}^{t}M+\gamma\mu'(s)\mu^{-1}(s)ds}\\
&=||x_{0}-x'_{0}||e^{M(t-t_{0})}
\left(\frac{\mu(t)}{\mu(t_{0})}\right)^{\gamma}.\\
||Y(t,t_{0},y_{0})-Y(t,t_{0},y'_{0})||&\leq ||y_{0}-y'_{0}|| e^{\int_{t_{0}}^{t}Mds}\\
&=||y_{0}-y'_{0}||e^{M(t-t_{0})}.
\end{align*}
\section{Proofs of main results}
{\bf Proof of Theorem {\ref{thm3.1}}.}
We are going to show that $H(t,\cdot)$ satisfies the four conditions of Definition \ref{def2.2}.\\
\indent For any fixed $t$, it follows from Lemmas \ref{lem4.7} and \ref{lem4.8} that $H(t, \cdot)$ is homeomorphism
and $G(t, \cdot) = H^{-1}(t, \cdot)$. Thus, Condition (i) is satisfied.\\
\indent From (\ref{4.26}) and Lemma \ref{lem4.2}, we derive $||H(t, x)-x|| = ||h(t,(t, x))|| \leq 2K\gamma\alpha^{-1}$. Note
$||H(t, x)||\to \infty$ as $|x| \to  \infty $, uniformly with respect to t. Thus, Condition (ii) is satisfied.\\
\indent From (\ref{4.27}) and Lemma \ref{lem4.3}, we derive $||G(t, y)-y|| = ||g(t,(t, y))|| \leq 2K\gamma\alpha^{-1}$. Note
$||G(t, y)||\to \infty$ as $|y| \to  \infty $, uniformly with respect to t. Thus, Condition (iii) is satisfied.\\
\indent From Lemmas \ref{lem4.5} and \ref{lem4.6}, we know that Condition (iv) is true.\\
\indent Hence, the system (\ref{2.1}) and its linear system (\ref{2.2}) are topologically conjugated.\\

{\bf Proof  of Theorem \ref{thm3.2}.} We prove this theorem in two steps.\\
{\bf step 1.} We show that there exist constants $p >0,0<q<1$, such that  $||H(t, x)-H(t, x)||\leq p||x-x'||^{q}$, if $||x-x'|| < 1$.

{\bf Proof.} From Lemma \ref{lem4.2}, it follows that
\begin{align*}
h(t,(t,\xi))=-\int_{-\infty}^{t}T(t,s)P(s)f(s,X(s,t,\xi))ds
+\int_{t}^{+\infty}T(t,s)Q(s)f(s,X(s,t,\xi))ds.
\end{align*}
Thus, we get
\begin{align*}
h(t,(t,\xi))-h(t,(t,\xi'))=&-\int_{-\infty}^{t}T(t,s)P(s)(f(s,X(s,t,\xi))
-f(s,X(s,t,\xi')))ds\\
&+\int_{t}^{+\infty}T(t,s)Q(s)(f(s,X(s,t,\xi))-f(s,X(s,t,\xi')))ds.\\
=&I_{1}+I_{2}.
\end{align*}
We suppose that $0<||\xi-\xi'||<1$. Taking $\tau=\frac{1}{M+\gamma}\ln\frac{1}{||\xi-\xi'||}$.\\
Now  divide $I_{1} , I_{2}$ into two parts:
\begin{align*}
I_{1}=\int_{-\infty}^{t-\tau}+\int_{t-\tau}^{t}=I_{11}+I_{12},\\
I_{2}=\int_{t}^{t+\tau}+\int_{t+\tau}^{+\infty}=I_{21}+I_{22}.
\end{align*}
Then, by using (\ref{2.7}) and (\ref{3.1}), we have
\begin{align*}
||I_{11}||&\leq\int_{-\infty}^{t-\tau}||T(t,s)P(s)||2\beta\mu'(s)\mu^{-1}(s)ds\\
&\leq\int_{-\infty}^{t-\tau}K\left(\frac{\mu(t)}{\mu(s)}\right)^{-\alpha}
2\beta\mu'(s)\mu^{-1}(s)ds\\
&\leq2K\beta\alpha^{-1}\left(\frac{\mu(t-\tau)}{\mu(t)}\right)^{\alpha}.
\end{align*}
Taking sufficiently large constant $\tilde{M}$, such that
\begin{align*} \tilde{M}\geq \max\{[\log_\frac{\mu(t-\tau)}{\mu(t)}{||\xi-\xi'||}]+1,[\log_\frac{\mu(t)}{\mu(t+\tau)}{||\xi-\xi'||}]+1,\alpha,M\}.
\end{align*}
Then, we get
\begin{align*}
\left(\frac{\mu(t-\tau)}{\mu(t)}\right)^{\alpha}\leq||\xi-\xi'||^{\frac{\alpha}{\tilde{M}}}.
\end{align*}
Thus,
\begin{align*}
||I_{11}||\leq2K\beta\alpha^{-1}||\xi-\xi'||^{\frac{\alpha}{\tilde{M}}},
\end{align*}
and
\begin{align*}
||I_{22}||&\leq\int_{t+\tau}^{+\infty}||T(t,s)Q(s)||2\beta\mu'(s)\mu^{-1}(s)ds\\
&\leq\int_{t+\tau}^{+\infty } K\left(\frac{\mu(s)}{\mu(t)}\right)^{-\alpha}2\beta\mu'(s)\mu^{-1}(s)ds\\
&\leq2K\beta\alpha^{-1}\left(\frac{\mu(t)}{\mu(t+\tau)}\right)^{\alpha}\\
&\leq2K\beta\alpha^{-1}||\xi-\xi'||^{\frac{\alpha}{\tilde{M}}}.
\end{align*}
From (\ref{2.7}), (\ref{3.1}) and Lemma  {\ref{lem4.10}}, we have
\begin{align*}
||I_{12}||
&\leq\int_{t-\tau}^{t}K\left(\frac{\mu(t)}{\mu(s)}\right)^{-\alpha}\gamma\mu'(s)\mu^{-1}(s)
||X(s,t,\xi)-X(s,t,\xi')||ds\\
&\leq\int_{t-\tau}^{t}K\left(\frac{\mu(t)}{\mu(s)}\right)^{-\alpha}\gamma\mu'(s)\mu^{-1}(s)
||\xi-\xi'||e^{M(s-t)}\left(\frac{\mu(s)}{\mu(t)}\right)^{\gamma}ds\\
&\leq K\gamma||\xi-\xi'||e^{M\tau}\int_{t-\tau}^{t}\left(\frac{\mu(t)}{\mu(s)}\right)^{- (\alpha+\gamma)}\mu'(s)\mu^{-1}(s)
ds\\
&\leq K\gamma||\xi-\xi'||e^{M\tau}(\alpha+\gamma)^{-1}\\
&\leq K\gamma(\alpha+\gamma)^{-1}||\xi-\xi'||^{\frac{\gamma}{M+\gamma}},
\end{align*}
and
\begin{align*}
||I_{21}||&\leq\int_{t}^{t+\tau}K\left(\frac{\mu(s)}{\mu(t)}\right)^{-\alpha}\gamma\mu'(s)\mu^{-1}(s)
||X(s,t,\xi)-X(s,t,\xi')||ds\\
&\leq\int_{t}^{t+\tau}K\left(\frac{\mu(s)}{\mu(t)}\right)^{-\alpha}\gamma\mu'(s)\mu^{-1}(s)
||\xi-\xi'||e^{M(s-t)}\left(\frac{\mu(s)}{\mu(t)}\right)^{\gamma}ds\\
&\leq K\gamma||\xi-\xi'||e^{M\tau}\int_{t}^{t+\tau}\left(\frac{\mu(s)}{\mu(t)}\right)^{-(\alpha-\gamma)}\mu'(s)\mu^{-1}(s)
ds\\
&\leq K\gamma||\xi-\xi'||e^{M\tau}(\alpha-\gamma)^{-1}\\
&\leq K\gamma(\alpha-\gamma)^{-1}||\xi-\xi'||^{\frac{\gamma}{M+\gamma}}.
\end{align*}
By the definition of $H(t,x)$, if $||x-x'||<1$,
\begin{align*}
||H(t,x)-H(t,x')||&\leq||x-x'||+||I_{11}||+||I_{12}||+||I_{21}||+||I_{22}||\\
&\leq(1+4K\beta\alpha^{-1}+K\gamma(\alpha+\gamma)^{-1}+K\gamma (\alpha-\gamma)^{-1})||x-x'||^{\min\{\frac{\alpha}{\tilde{M}},\frac{\gamma}{M+\gamma}\}}\\
&\equiv p||x-x'||^{q},
\end{align*}
where $p=1+4K\beta\alpha^{-1}+K\gamma(\alpha+\gamma)^{-1}+K\gamma (\alpha-\gamma)^{-1})$ , $q=\min\{\frac{\alpha}{\tilde{M}},\frac{\gamma}{M+\gamma}\}.$

{\bf Step 2.} We show that there exist constants $p'>0, 0<q'<1$, such that $||G(t, y)-G(t, y')||\leq p'||y-y'||^{q'}$, if $||y-y'|| < 1$.

{\bf Proof.} From Lemma \ref{lem4.3}, we know that $g(t,(\tau, \xi))$ is a fixed point of  map $\mathscr{F}$ :
\begin{align*}
\mathscr{F} z=\int_{-\infty}^{t}T(t,s)P(s)f(s,Y(s,\tau,\xi)+z)ds
-\int_{t}^{+\infty}T(t,s)Q(s)f(s,Y(s,\tau,\xi)+z)ds.
\end{align*}
Let $g_{0}(t,(\tau,\xi))\equiv0$ and recursively define:
\begin{align*}
g_{m+1}(t,(\tau,\xi))=&\int_{-\infty}^{t}T(t,s)P(s)f(s,Y(s,\tau,\xi)+g_{m}(s,(\tau,\xi)))ds\\
&-\int_{t}^{+\infty}T(t,s)Q(s)f(s,Y(s,\tau,\xi)+g_{m}(s,(\tau,\xi)))ds.
\end{align*}
Firstly, we prove that
\begin{align*}
g_{m}(t,(\tau,\xi)) \to g(t,(\tau,\xi)),\  \mathrm{as}\  m \to +\infty,
\end{align*}
uniformly with respect to $t,\tau,\xi$.
\begin{align*}
||g_{1}(t,(\tau,\xi))-g_{0}(t,(\tau,\xi))||\leq&\int_{-\infty}^{t}||T(t,s)P(s)
f(s,Y(s,\tau,\xi)+g_{0}(s,(\tau,\xi))||ds\\
&-\int_{t}^{+\infty}||T(t,s)Q(s)f(s,Y(s,\tau,\xi)+g_{0}(s,(\tau,\xi))||ds\\
\leq&\int_{-\infty}^{t}K\left(\frac{\mu(t)}{\mu(s)}\right)^{-\alpha}\beta\mu'(s)
\mu^{-1}(s)ds\\
&+\int_{t}^{+\infty}K\left(\frac{\mu(s)}{\mu(t)}\right)^{-\alpha}\beta\mu'(s)
\mu^{-1}(s)|ds\\
\leq&4K\beta\alpha^{-1},
\end{align*}
\begin{align*}
||g_{2}(t,(\tau,\xi))-g_{1}(t,(\tau,\xi))||\leq&\int_{-\infty}^{t}||T(t,s)P(s)
(f(s,Y(s,\tau,\xi)+g_{1}(s,(\tau,\xi))))\\
&-f(s,Y(s,\tau,\xi)+g_{0}(s,(\tau,\xi))))||ds\\
&-\int_{t}^{+\infty}||T(t,s)Q(s)f(s,Y(s,\tau,\xi)+g_{1}(s,(\tau,\xi))\\
&-f(s,Y(s,\tau,\xi)+g_{0}(s,(\tau,\xi))))||ds\\
\leq&\int_{-\infty}^{t}K\left(\frac{\mu(t)}{\mu(s)}\right)^{-\alpha}\gamma\mu'(s)
\mu^{-1}(s)||g_{1}-g_{0}||ds\\
&+\int_{t}^{+\infty}K\left(\frac{\mu(s)}{\mu(t)}\right)^{-\alpha}\gamma\mu'(s)
\mu^{-1}(s)||g_{1}-g_{0}||ds\\
\leq&(4K\gamma\alpha^{-1})(4K\beta\alpha^{-1}).
\end{align*}
Similarly, we get
\begin{align*}
||g_{m+1}(t,(\tau,\xi))-g_{m}(t,(\tau,\xi))||\leq&\int_{-\infty}^{t}||T(t,s)P(s)
(f(s,Y(s,\tau,\xi)+g_{m})\\
&-f(s,Y(s,\tau,\xi)+g_{m-1})||ds\\
&-\int_{t}^{+\infty}||T(t,s)Q(s)f(s,Y(s,\tau,\xi)+g_{m})\\
&-f(s,Y(s,\tau,\xi)+g_{m-1})||ds\\
\leq&\int_{-\infty}^{t}K\left(\frac{\mu(t)}{\mu(s)}\right)^{-\alpha}\gamma\mu'(s)
\mu^{-1}(s)||g_{m}-g_{m-1}||ds\\
&+\int_{t}^{+\infty}K\left(\frac{\mu(s)}{\mu(t)}\right)^{-\alpha}\gamma\mu'(s)
\mu^{-1}(s)||g_{m}-g_{m-1}||ds\\
\leq&(4K\gamma\alpha^{-1})^{m}(4K\beta\alpha^{-1}).
\end{align*}
Since $4K\gamma\alpha^{-1}<1$, $g_{m}(t,(\tau,\xi))$ is uniformly converged, as $m\to +\infty$. Let
$$\lim_{m\to +\infty}g_{m}(t,(\tau,\xi))=\bar{g}(t,(\tau,\xi)).$$
Secondly, we prove $\bar{g}(t,(\tau,\xi))=g(t,(\tau,\xi))$.
\begin{align*}
||g_{0}(t,\tau,\xi)-g(t,\tau,\xi)||\leq&\int_{-\infty}^{t}||T(t,s)P(s)
f(s,Y(s,\tau,\xi)+g)||ds\\
&-\int_{t}^{+\infty}||T(t,s)Q(s)f(s,Y(s,\tau,\xi)+g)ds||\\
\leq&\int_{-\infty}^{t}K\left(\frac{\mu(t)}{\mu(s)}\right)^{-\alpha}\beta\mu'(s)
\mu^{-1}(s)ds\\
&+\int_{t}^{+\infty}K\left(\frac{\mu(s)}{\mu(t)}\right)^{-\alpha}\beta\mu'(s)
\mu^{-1}(s)ds\\
\leq&4K\beta\alpha^{-1},
\end{align*}
\begin{align*}
||g_{1}(t,(\tau,\xi))-g(t,(\tau,\xi))||\leq&\int_{-\infty}^{t}||T(t,s)P(s)
(f(s,Y(s,\tau,\xi)+g_{0}(s,(\tau,\xi)))\\
&-f(s,Y(s,\tau,\xi)+g(s,(\tau,\xi))))||ds\\
&-\int_{t}^{+\infty}||T(t,s)Q(s)f(s,Y(s,\tau,\xi)+g_{0}(s,(\tau,\xi))\\
&-f(s,Y(s,\tau,\xi)+g(s,(\tau,\xi))))||ds\\
\leq&\int_{-\infty}^{t}K\left(\frac{\mu(t)}{\mu(s)}\right)^{-\alpha}\gamma\mu'(s)
\mu^{-1}(s)||g_{0}-g||ds\\
&+\int_{t}^{+\infty}K\left(\frac{\mu(s)}{\mu(t)}\right)^{-\alpha}\gamma\mu'(s)
\mu^{-1}(s)||g_{0}-g||ds\\
\leq&(4K\gamma\alpha^{-1})(4K\beta\alpha^{-1}).
\end{align*}
Similarly, we get
\begin{align*}
||g_{m}(t,(\tau,\xi))-g(t,(\tau,\xi))||\leq&\int_{-\infty}^{t}||T(t,s)P(s)
(f(s,Y(s,\tau,\xi)+g_{m-1})\\
&-f(s,Y(s,\tau,\xi)+g))||ds\\
&-\int_{t}^{+\infty}||T(t,s)Q(s)(f(s,Y(s,\tau,\xi)+g_{m-1})\\
&-f(s,Y(s,\tau,\xi)+g))||ds\\
\leq&\int_{-\infty}^{t}K\left(\frac{\mu(t)}{\mu(s)}\right)^{-\alpha}\gamma\mu'(s)
\mu^{-1}(s)||g-g_{m-1}||ds\\
&+\int_{t}^{+\infty}K\left(\frac{\mu(s)}{\mu(t)}\right)^{-\alpha}\gamma\mu'(s)
\mu^{-1}(s)||g-g_{m-1}||ds\\
\leq&(4K\gamma\alpha^{-1})^{m}(4K\beta\alpha^{-1}).
\end{align*}
Since $4K\gamma\alpha^{-1}<1$,   $g_{m}(t,(\tau,\xi))\to g(t,(\tau,\xi))$ , as $m\to +\infty$.
From the uniqueness of the limit, we know that $\bar{g}(t,(\tau,\xi))=g(t,(\tau,\xi))$. Thus,
$$\lim_{m\to+\infty}g_{m}(t,\tau,\xi)=g(t,\tau,\xi).$$
Next, we prove $g_{m}(t,(\tau,\xi))=g_{m}(t,(t,Y(t,\tau,\xi)))$. Since
$$Y(t,(\tau,\xi))=Y(t,t,Y(t,\tau,\xi)).$$
$$g_{1}(t,(\tau,\xi))=\int_{-\infty}^{t}T(t,s)P(s)f(s,Y(s,\tau,\xi))ds-
\int_{t}^{+\infty}T(t,s)Q(s)f(s,Y(s,\tau,\xi))ds.$$
\begin{equation*}
\begin{split}
g_{1}(t,(t,Y(t,\tau,\xi)))=&\int_{-\infty}^{t}T(t,s)P(s)f(s,Y(s,s,Y(s,\tau,\xi)))ds\\
&-\int_{t}^{+\infty}T(t,s)Q(s)f(s,Y(s,s,Y(s,\tau,\xi)))ds\\
=&\int_{-\infty}^{t}T(t,s)P(s)f(s,Y(s,\tau,\xi))ds
-\int_{t}^{+\infty}T(t,s)Q(s)f(s,Y(s,\tau,\xi))ds.
\end{split}
\end{equation*}
Thus, $g_{1}(t,(\tau,\xi))=g_{1}(t,(t,Y(t,\tau,\xi))).$\\
Similarly, we get
\begin{align*}
g_{m+1}(t,(\tau,\xi))=&\int_{-\infty}^{t}T(t,s)P(s)f(s,Y(s,\tau,\xi)+g_{m}(s,\tau,\xi))ds\\
&-\int_{t}^{+\infty}T(t,s)Q(s)f(s,Y(s,\tau,\xi)+g_{m}(s,\tau,\xi))ds,\\
g_{m+1}(t,(t,Y(t,\tau,\xi)))=&\int_{-\infty}^{t}T(t,s)P(s)f(s,Y(s,s,Y(s,\tau,\xi))+g_{m}(s,s,Y(s,\tau,\xi)))ds
\\&-\int_{t}^{+\infty}T(t,s)Q(s)f(s,Y(s,s,Y(s,\tau,\xi))+g_{m}(s,s,Y(s,\tau,\xi)))ds\\
=&\int_{-\infty}^{t}T(t,s)P(s)f(s,Y(s,\tau,\xi)+g_{m}(s,\tau,\xi))ds\\
&-\int_{t}^{+\infty}T(t,s)Q(s)f(s,Y(s,\tau,\xi))ds.
\end{align*}
Thus, $g_{m+1}(t,(\tau,\xi))=g_{m+1}(t,(t,Y(t,\tau,\xi))).$\\
From the mathematical induction, we get
$$g_{m}(t,(\tau,\xi))=g_{m}(t,(t,Y(t,\tau,\xi))).$$
Similarly, we take sufficiently large constant $\tilde{M_{1}}$ and $\tilde{\tau}=\frac{1}{M+\gamma}\ln\frac{1}{||\eta-\eta'||}$, such that
\begin{align*}
\tilde{M_{1}}\geq \max\{[\log_\frac{\mu(t-\tilde{\tau})}{\mu(t)}{||\eta-\eta'||}]+1,[\log_\frac{\mu(t)}{\mu(t+\tilde{\tau})}{||\eta-\eta'||}]+1,\alpha,M\}.
\end{align*}
Finally, taking constants $\lambda>0$ is sufficiently large, and $q'>0$ is sufficiently small, satisfing  $\lambda>\frac{4K\beta\alpha^{-1}+2K\gamma\alpha^{-1}}{1-2K\gamma\alpha^{-1}e}$, $q'\leq \min\{\frac{\alpha}{\tilde{M_{1}}},\frac{\gamma}{M+\gamma},\frac{\alpha}{M\tilde{\tau}}\ln\frac{\mu(t+\tilde{\tau})}{\mu(t)}\}$, where $\ln\frac{\mu(t+\tilde{\tau})}{\mu(t)}<\frac{1}{\alpha}$.\\
Furthermore, we prove  if $ ||\eta-\eta'|| < 1$, for any non-negative integer $m$, we have
\begin{align}\label{foruma4.14}
||g_{m}(t,(t,\eta))-g_{m}(t,(t,\eta'))||\leq \lambda||\eta-\eta'||^{q'}.
\end{align}
When $m=0$, (\ref{foruma4.14}) obviously holds. Now we induce that hypothesis (\ref{foruma4.14}) holds. From (\ref{foruma4.14}), we get
\begin{align*}
g_{m+1}(t,(t,\eta))-g_{m+1}(t,(t,\eta'))=&
\int_{-\infty}^{t}T(t,s)P(s)[f(s,Y(s,t,\eta)+g_{m}(s,t,\eta))\\
&-f(s,Y(s,t,\eta')+g_{m}(s,t,\eta'))]ds\\
&-\int_{t}^{+\infty}T(t,s)Q(s)[f(s,Y(s,t,\eta)+g_{m}(s,t,\eta))\\
&-f(s,Y(s,t,\eta')+g_{m}(s,t,\eta'))]ds\\
=&J_{1}+J_{2}.
\end{align*}
Now  divide $J_{1} , J_{2}$ into two parts:
\begin{align*}
J_{1}=\int_{-\infty}^{t-\tilde{\tau}}+\int_{t-\tilde{\tau}}^{t}=J_{11}+J_{12},\\
J_{2}=\int_{t}^{t+\tilde{\tau}}+\int_{t+\tilde{\tau}}^{+\infty}=J_{21}+J_{22}.
\end{align*}
Then, by using (\ref{2.7}) and (\ref{3.1}), we have
\begin{align*}
||J_{11}||&\leq\int_{-\infty}^{t-\tilde{\tau}}||T(t,s)P(s)||2\beta\mu'(s)\mu^{-1}(s)ds\\
&\leq\int_{-\infty}^{t-\tilde{\tau}}K\left(\frac{\mu(t)}{\mu(s)}\right)^{-\alpha}
2\beta\mu'(s)\mu^{-1}(s)ds\\
&\leq2K\beta\alpha^{-1}\left(\frac{\mu(t-\tilde{\tau})}{\mu(t)}\right)^{\alpha}\\
&\leq2K\beta\alpha^{-1}||\eta-\eta'||^{\frac{\alpha}{\tilde{M_{1}}}},
\end{align*}
and
\begin{align*}
||J_{22}||&\leq\int_{t+\tilde{\tau}}^{+\infty}||T(t,s)Q(s)||2\beta\mu'(s)\mu^{-1}(s)ds\\
&\leq\int_{t+\tilde{\tau}}^{+\infty } K\left(\frac{\mu(s)}{\mu(t)}\right)^{-\alpha}2\beta\mu'(s)\mu^{-1}(s)ds\\
&\leq2K\beta\alpha^{-1}\left(\frac{\mu(t)}{\mu(t+\tilde{\tau})}\right)^{\alpha}\\
&\leq2K\beta\alpha^{-1}||\eta-\eta'||^{\frac{\alpha}{\tilde{M_{1}}}}.
\end{align*}
Furthermore, it follows from (\ref{2.7}) , (\ref{3.1})  and lemma {\ref{lem4.10}} that
\begin{align*}
||J_{12}||
&\leq\int_{t-\tilde{\tau}}^{t}K\left(\frac{\mu(t)}{\mu(s)}\right)^{-\alpha}\gamma\mu'(s)\mu^{-1}(s)
(||\eta-\eta'||e^{M(s-t)}+\lambda||\eta-\eta'||^{q'}e^{Mq'(s-t)})ds\\
&\leq\int_{t-\tilde{\tau}}^{t}K\left(\frac{\mu(t)}{\mu(s)}\right)^{-\alpha}\gamma\mu'(s)\mu^{-1}(s)
(||\eta-\eta'||e^{M\tilde{\tau}}+\lambda||\eta-\eta'||^{q'}e^{Mq'\tilde{\tau}})ds\\
&\leq K\gamma\alpha^{-1}||\eta-\eta'||e^{M\tilde{\tau}}+K\gamma\lambda||\eta-\eta'||^{q'}\left(\frac{\mu(t)}{\mu(t-\tilde{\tau})}\right)^{\alpha}\ln\frac{\mu(t)}{\mu(t-\tilde{\tau})}\\
&\leq K\gamma\alpha^{-1}||\eta-\eta'||^{\frac{\gamma}{M+\gamma}}+K\gamma\lambda\alpha^{-1}e||\eta-\eta'||^{q'}\\
&\leq (K\gamma\alpha^{-1}+K\gamma\lambda\alpha^{-1}e)||\eta-\eta'||^{q'},
\end{align*}
and
\begin{align*}
||J_{21}||
&\leq\int_{t}^{t+\tilde{\tau}}K\left(\frac{\mu(s)}{\mu(t)}\right)^{-\alpha}\gamma\mu'(s)\mu^{-1}(s)
(||\eta-\eta'||e^{M(s-t)}+\lambda||\eta-\eta'||^{q'}e^{Mq'(s-t)})ds\\
&\leq\int_{t}^{t+\tilde{\tau}}K\left(\frac{\mu(s)}{\mu(t)}\right)^{-\alpha}\gamma\mu'(s)\mu^{-1}(s)
(||\eta-\eta'||e^{M\tilde{\tau}}+\lambda||\eta-\eta'||^{q'}e^{Mq'\tilde{\tau}})ds\\
&\leq K\gamma\alpha^{-1}||\eta-\eta'||e^{M\tilde{\tau}}+K\gamma\lambda||\eta-\eta'||^{q'}\left(\frac{\mu(t+\tilde{\tau})}{\mu(t)}\right)^{\alpha}\ln\frac{\mu(t+\tilde{\tau})}{\mu(t)}\\
&\leq K\gamma\alpha^{-1}||\eta-\eta'||^{\frac{\gamma}{M+\gamma}}+K\gamma\lambda\alpha^{-1}e||\eta-\eta'||^{q'}\\
&\leq (K\gamma\alpha^{-1}+K\gamma\lambda\alpha^{-1}e)||\eta-\eta'||^{q'}.
\end{align*}
Thus,
\begin{align*}
||g_{m+1}(t,(t,\xi))-g_{m+1}(t,(t,\xi'))||&\leq||J_{11}||+||J_{12}||+||J_{21}||+||J_{22}||\\
&\leq4K\beta\alpha^{1}||\eta-\eta'||^{\frac{\alpha}{\tilde{M_{1}}}}+2(K\gamma\alpha^{-1}+K\gamma\lambda\alpha^{-1}e)||\eta-\eta'||^{q'}\\
&\leq\left(4K\beta\alpha^{-1}+2K\gamma\alpha^{-1}+2K\gamma\alpha^{-1}\lambda e\right)||\eta-\eta'||^{q'}\\
&\leq\lambda||\eta-\eta'||^{q'}.
\end{align*}
This means that for any non-negative integer $m$, there is
\begin{align*}
||g(t,(t,\eta))-g(t,(t,\eta'))||\leq\lambda||\eta-\eta'||^{q'}.
\end{align*}
From the definition of $G(t,y)$, if $0<||y-y'||<1$, we have
\begin{align*}
||G(t,y)-G(t,y')||\leq||y-y'||+\lambda||y-y'||^{q'}\leq(1+\lambda)||y-y'||^{q'}.
\end{align*}
This completes the proof of Theorem \ref{thm3.2}.
\section{An example}
Consider the differential equation in $\mathbb{R}^2$,
\begin{equation}\label{exam5.1}
\begin{split}
\begin{cases}
x'=-\eta_{1}\left(\frac{\mu'(t)}{\mu(t)}\right)x+\varepsilon \sin(x+t),\\
y'=\eta_{2}\left(\frac{\mu'(t)}{\mu(t)}\right)y+\varepsilon \cos(x+t),
\end{cases}
\end{split}
\end{equation}
for $t\in \mathbb{R}$, where $\eta_{1}, \eta_{2}, \varepsilon$ are positive constants and $\mu(t)=\frac{2}{\pi}e^{t}(\frac{\pi}{2}+\arctan t)$. Denoting
\begin{equation*}
z=
\left(
\begin{array}{cc}
x\\
y\\
\end{array}
\right),\
A(t)=
\left(
\begin{array}{cc}
-\eta_{1}\left(\frac{\mu(t)}{\mu(s)}\right)&0\\
0&\eta_{2}\left(\frac{\mu'(t)}{\mu(t)}\right)\\
\end{array}
\right),\
f(t,z)=
\left(
\begin{array}{cc}
\varepsilon \sin(x+t)\\
\varepsilon \cos(x+t)\\
\end{array}
\right).
\end{equation*}
From example \ref{exam2.2}, we know that the linear part of equation (\ref{exam5.1}) admits an algebraic dichotomy.
Then,
\begin{align*}
\mu'(t)\mu^{-1}(t)=1+\frac{1}{(\frac{\pi}{2}+\arctan t)(1+t^{2})}.
\end{align*}
We can easily get $1<\mu'(t)\mu^{-1}(t)<1+\frac{2}{\pi}$.
Then, we have
\begin{align*}
||f(t,z)||\leq\varepsilon||\sin(x+t)+\cos(x+t)||\leq2\varepsilon\leq2\varepsilon\mu'(t)\mu^{-1}(t),\\
||f(t,z_{1})-f(t,z_{2})||\leq2\varepsilon||z_{1}-z_{2}||\leq2\varepsilon\mu'(t)\mu^{-1}(t).
\end{align*}
Hence, equation (\ref{exam5.1}) satisfies the condition of Theorem \ref{thm3.1} if $0<\varepsilon<\frac{\eta_{3}}{8}$, where $\eta_{3}=\max\{\eta_{1},\eta_{2}\}$. Therefore, Eq. (\ref{exam5.1}) is topologically conjugated to its linear part.

\section{Data Availability Statement}
%\hskip\parindent
My manuscript has no associated data. It is pure mathematics.

\section{Conflict of Interest}
%\hskip\parindent
The authors declare that they have no conflict of interest.

\subsection{Authors' contributions}

All authors contributed equally to the writing of this paper. All authors read and approved the final manuscript.

\end{document}